\documentclass[12pt]{article}

\usepackage{amsfonts}
\usepackage{amssymb}
\def\Bas{\mathop{\rm Bas}}

\def\charPol{\mathop{\rm charPol}}
\def\Cir{\mathop{\rm Cir}}
\def\C1G{\mathop{\rm C1G}}
\def\codim{\mathop{\rm codim}}

\def\Dg{\mathop{\rm Dg}}
\def\diag{\mathop{\rm diag}}
\def\dim{\mathop{\rm dim}}
\def\dimSO{\mathop{\rm dimSO}}
\def\discr{\mathop{\rm discr}}
\def\DO{\mathop{\rm DO}}

\def\elSym{\mathop{\rm elSym}}
\def\ev{\mathop{\rm ev}}

\def\Orbit{\mathop{\rm Orbit}}
\def\Rt{\mathop{\rm Rt}}
\def\RtA{\mathop{\rm RtA}}
\def\Scal{\mathop{\rm Scal}}
\def\scp{\mathop{\rm scp}}
\def\sgn{\mathop{\rm sgn}}
\def\so{\mathop{\rm so}}
\def\SO{\mathop{\rm SO}}
\def\spe{\mathop{\rm spe}}
\def\spec{\mathop{\rm spec}}
\def\sQuad{\mathop{\rm sQuad}}
\def\St{\mathop{\rm St}}
\def\Sym{\mathop{\rm Sym}}
\def\trace{\mathop{\rm trace}}
\begin{document}
\title{Geometry of the Discriminant Surface\\ for Quadratic Forms}
\author{S.~D.~Mechveliani\\
        Program Systems Institute of Russian Academy of Science,\\
        Pereslavl-Zalessky, Russia. \ \ 
        e-mail: \ mechvel@botik.ru
       }
\date{October 4, 2011.}
\maketitle
\begin{abstract}
We investigate the manifold \ $\cal{M}$ \ of (real) quadratic forms in \ 
$n > 1$ \ variables having a multiple eigenvalue. \ \
In addition to known facts, we prove that\linebreak 
1) $\cal{M}$ is irreducible, \ \ 
2) in the case of $n = 3$, \ scalar matrices and only them are singular
   points on $\cal{M}$. 

For $n = 3$, \ $\cal{M}$ is also described as the straight cylinder over 
$\cal{M}$$_0$, where $\cal{M}$$_0$ is the cone over the orbit of the 
diagonal matrix $\diag(1,1,-2)$ by the orthogonal changes of coordinates.
We analyze certain properties of this orbit, which occurs a diffeomorphic 
image of the projective plane. \ Bibliography: 9 items.
\end{abstract}

                  singular point.
\section{ Introduction }
\label{sec-intr}


This paper is an improvement and abbreviation of the manuscript \cite{Me}
version of 2010.

\noindent {\bf Definition}
 
1) \ $\Sym(n) = \Sym(n, \mathbb{R})$ \ is the space of symmetric matrices
   of size $n$ over the field \ $\mathbb{R}$ \ of real numbers 
   ($\dim \Sym(n) = n (n+1)/2$).

  Symmetric matrices \ we also call \ {\it quadratic forms} \ or \ 
  {\it forms} \ \ --- presuming that a quadratic form (a bilinear symmetric 
  function) is written in coordinates as a symmetric matrix.

2) \ $\charPol(X)(\lambda)$ \ \ 
                          is the characteristic polynomial of a matrix $X$,   

  $\discr(X)$ \ is discriminant of \ $\charPol(X)$.

3) \ $\cal{M}$ $ = \{X \in \Sym(n) | \ \discr(X) = 0\}$ \ \ 
   is the {\it surface} in $\Sym(n)$ defined by discriminant. \
   In other words: \ $\cal{M}$ is the surface of all the forms in 
   $\Sym(n)$ that have any multiple eigenvalue. 

4) By the word ``surface'' \ we mean an \ open subset in a smooth manifold 
   or in an \ affine algebraic variety \ in a finite-dimensional Euclidean
   space \ (a surface may have singular points).  

5) The word \ ``fibration'' \ \ means in this paper \ a fibration over a plane 
   considered \ a) modulo isometry, \ \ b) modulo diffeomorphism.

6) \ $\cal{MH}$ \ is the set of the forms in $\Sym(n)$ \ having exactly \  
   $n-1$ \ different eigenvalues.
\smallskip

In the book \cite{Ar1}, in Appendix 10, \ the open area $\Sym_+(n)$ of the
positive definite forms in $\Sym(n)$ is represented as the area of ellipsoids.
And it is shown that (in our denotations) \ \ \ 
a) $\codim \cal{M}$ $= 2$, \ \ \ b) $\cal{MH}$ is a smooth manifold.

Also it simply follows from the theorem and lemma in this paragraph of 
\cite{Ar1} that $\cal{MH}$ is an open in $\cal{M}$ and dense in $\cal{M}$ 
set.
\smallskip

\noindent \underbar{Aim of this paper:} \ \ to answer the questions about 

{\bf (1)} algebraic irreducibility of $\cal{M}$, 

{\bf (2)} location of its singular points, its global structure in 
          the case of \ $n = 3$,

{\bf (3)} algebraic equations for the orbit of a form in $\cal{M}$ by 
          the orthogonal coordinate changes, and of its position 
          relatively to the plane of diagonal matrices. 

As for \ $n = 1$ \ the set $\cal{M}$ is empty, the problem is considered
only for \ $n > 1$.

It is known that each matrix in $\Sym(n, \mathbb{R})$ has \ $n$ \ real 
eigenvalues, some of which may coincide (see below Theorem DO).

\noindent {\bf Definition.} \ 
We call the \ {\it spectrum} \ or \ the \ {\it multiset of eigenvalues} \ 
the set of eigenvalues for a matrix in $\Sym(n)$ --- together with their 
multiplicities, \ and denote it \ \ $\spec$.

\subsection{About applications}

In Appendix 10 of \cite{Ar1} 
(``Multiplicities of characteristic frequencies and ellipsoids
depending on parameters''), \ it is described the mechanical problem of 
adjusting of a rigid body by moving of an adjustable mass along an arc 
rigidly attached to the body. The aim is to make the {\it ellipsoid of inertia} 
of the body to become an {\it ellipsoid of revolution} (=== of rotation).
And it turns out that a single adjustable mass (single parameter) is not 
sufficient in the generic configuration.
The reason for this is that codimension of $\mathbb{M}$ occurs not 1 
(as one could expect) but 2.

In the section 6 of \cite{Ar2} it is written about the relevance of this 
problem to the von-Neumann--Wigner theorem about intersection of electronic 
levels. 

In our understanding, from the point of view of mechanics and physics, a 
quadratic form $X(p)$ depending on a parameter $p$, represents some physical 
characteristic of an object, where the eigenvalues are the oscillation 
frequencies of some parts of an object (construction). 
Hence the discriminant surface represents a geometrical condition of the 
resonance effect.

Generally: the discriminant surface is an important object, and the more is 
known about it the better.

\subsection{Acknowledgements} 
I am grateful to \ {\it A.~G.~Khovanskii} \ for the reference to Arnold's
book (due to which it has been cancelled the large paper part which re-proved 
the results from this book). \ 
I am grateful to \ {\it N.~V.~Ilyushechkin, S.~Yu.~Orevkov, Yu.~L.~Sachkov} \ 
for certain bits of oral discussion.

\subsection{Case $n = 2$}

For the quadratic form \ 
{\footnotesize $$X = \pmatrix{x & y \cr
                              y & z}$$
}\noindent 
we have \
$\charPol(X)(\lambda) = \ (x - \lambda)(z - \lambda) - y^2 = \ \               
\lambda^2 - (x+z) \lambda + x z - y^2$.
\\
Hence discriminant of this form (as a polynomial in the elements of $X$) is
$$(x+z)^2 - 4 (x z - y^2) = \                                        
  x^2 + 2 x z + z^2 - 4 x z + 4 y^2 = \ \ (x - z)^2 + (2 y)^2.
$$
\noindent
It follows from this decomposition into a sum of squares that \ $\cal{M}$ \ 
is defined by the equations \ \ $\{x = z, \ y = 0\}$ \ \ and is the
{\it line of scalar matrices}: \ an irreducible surface of codimension 2 in 
$\Sym(2)$. \ \ All the points are regular in this case. \  
So, for the case of $n = 2$ \ the problem is solved. \ 
{\it In the sequel, it is considered only the case of} \ $n > 2$.

\subsection{About equations for $\cal{M}$}

For $n = 3$, \ $\discr(X)$ is a homogeneous form of degree 6 in 6 variables
$x_{i,j}$, consisting of 123 non-zero monomials. How to derive from this
expression a knowledge about dimension of $\cal{M}$, its irreducible 
components, singular points? A naive consideration\linebreak  
``a single equation, therefore \ $\codim({\cal M}) = 1$'' \ is erroneous.

{\it The quadratic form discriminant is a sum of squares of several real
polynomials}. This has been independently proved by several authors,
starting from Borhardt in 1846 up to N.~V.~Ilyushechkin \cite{Il} in 1985.
Further, in the paper \cite{Il2} there are derived for the case of \ $n = 3$ \ 
certain relations between the members of this decomposition, and by this, it 
derives for the case of \ $n = 3$ \ certain four cubic homogeneous forms 
defining the discriminant surface. 

And this our paper completes the above results with various important
knowledge about the discriminant surface.

\subsection{Our approach}
\label{sec-intr:approach}

This is the classical method of a linear Lie group of symmetries.
All symmetric matrices are produced by conjugating diagonal matrices with 
operators from $\SO(n)$. A displace in the plane of diagonal matrices is 
orthogonal (in a certain Euclidean metric) to a displace along the conjugation 
orbit. This provides a smooth parameterization for $\cal{M}$ --- except certain 
singular points. 
It remains to describe the orbit for a diagonal matrix, with considering its
spectrum width: \ $n-1$ \ or less than \ $n-1$.
\smallskip

\noindent {\bf Definition}

* The metric on the surface $\cal{M}$ is taken from the space $\Sym(n)$. 

* The \ spectrum width \ for a form from $\Sym(n)$ \ is the number of its 
  different eigenvalues.

  For a set $\cal{X}$ in $\Sym(n)$, \  we call \ 
  {\it a form of maximal spectrum} in $\cal{X}$ \ 
  any form from $\cal{X}$ having maximal in $\cal{X}$ spectrum width, \
  and all other forms are called the forms of \ {\it narrowed spectrum}.

* \ $\cal{MH}$ \ denotes the set of forms of maximal spectrum in $\cal{M}$ 
  --- that is of the spectrum width \ $n-1$.

\subsection{Results}
\label{sec-intr:result}

* It is proved that $\cal{M}$ is an irreducible algebraic surface \ 
  (Section \ref{sec-prime}).

* There are presented algebraic equations for the orbit of a form in 
  $\cal{MH}$ under the orthogonal coordinate changes 
  (Section \ref{sec-DSForm:orbit}).

  It is proved that this orbit intersects orthogonally the plane of diagonal
  matrices.

* For \ $n = 3$ \ \ it is proved that the singular points on $\cal{M}$ are 
  scalar matrices and only them \ (Section \ref{sec-DSForm:structure3}).

* For \ $n = 3$ \ \ we provide a more definite description for $\cal{M}$: \
  the straight cylinder over $\cal{M}$$_0$, where $\cal{M}$$_0$ is the cone 
  over a diffeomorphic image of the projective plane $\mathbb{R}P(2)$ \
  (Lemma Orbit1, Sections \ref{sec-DSForm:orbit3}, 
   \ref{sec-DSForm:structure3}). 

* For \ $n = 3$ \ there are analyzed certain additional properties of the 
  orbit by orthogonal changes for a form in $\cal{MH}$: its certain 
  circumference family having a common point, and the minimal dimension of 
  an embracing plane (Section \ref{sec-DSForm:orbit3}).
\smallskip

\noindent {\bf The next questions can be as follows.}

1) For \ $n > 3$, \ is it true that each matrix in $\cal{M}$ having narrowed 
   spectrum is a singular point?

2) The structure of $\cal{M}$ depends mainly on the structure of an orbit 
 $\Orbit(A,n)$ for $A \in \cal{M}$. Similarly as for $n = 3$, 
 it would have sense to find for $n = 4, 5$ \ some global description for
 $\Orbit(A,n)$, more definite than the one given in the sequel.
\smallskip

\noindent {\bf Agreement.} \ \
1. To reduce the main text volume, some routine proofs are moved to 
   Appendix, and the main part contains references to Appendix.
 And\\
 {\it the proofs in Appendix do \ not \ use the statements being proved 
      in the main part of the paper}.

2. The labels of theorems and lemmata are unique, their scope is the whole
   text, including Appendix.


\newpage
\section{Results from Arnold's book}
\label{sec-Ar}

\noindent {\bf Theorem (Ar).}

(1) $\cal{M}$ is a finite union of smooth manifolds of codimension 2 and
    higher in the space $\Sym(n)$.

(2) $\cal{MH}$ is a smooth manifold of codimension 2 in $\Sym(n)$.

(3) $\cal{MH}$ is a set open in $\cal{M}$ and dense in $\cal{M}$.

\noindent {\bf Proof.}\\
This theorem immediately follows from the material in Appendix 10 of
the book \cite{Ar1}. The points (1) and (2) is our reformulation, for in 
this paragraph of \cite{Ar1} there are considered only positive definite 
forms in $\Sym(n)$.

This theorem is derived from \cite{Ar1} as follows.

\noindent {\bf Definition.}

1) \ $\Sym_+(n)$ \ is the set of positive definite forms in $\Sym(n)$. 
   It is defined by the condition of positiveness of the eigenvalues and is
   an open area in $\Sym(n)$.

2) \ $\cal{M}_+ = \cal{M} \cap$ $\Sym_+(n), \ \ \               
      \cal{MH}_+ = \cal{MH} \cap$ $\Sym_+(n)$.
\smallskip

Each of the manifolds \ $\Sym_+(n), \ \cal{M}_+, \ \cal{MH}_+$ \ is 
isometrically mapped on itself by the coordinate changes from $\SO(n)$ \
(see below Lemma S). \ Further, adding of a scalar matrix 
(applying a {\it scalar shift}) \ isometrically maps each of the manifolds \ 
$\Sym(n), \cal{M}, \cal{MH}$ \ on itself (see below Lemma Orbit1). 
And for each form in $\Sym(n)$ some its neighbourhood is mapped by some 
scalar shift into $\Sym(n)_+$. \ 
Therefore, the local metrical properties of the surface $\cal{M}_+$ are the 
same as of $\cal{M}$. This way (by a scalar shift), the theorem in \cite{Ar1} 
for ellipsoids is generalized to the case of the surface $\cal{M}$ in 
$\Sym(n)$.

In the book \cite{Ar1}, in the Section (A) of Appendix 10, it is described
a diffeomorphism from $\cal{M}_+$ onto the manifold of ellipsoids in the
space \ $V = \mathbb{R}^n$ \ having center in 0:\linebreak
        $A \ \mapsto \ \{v \in V | \ v A v* = 1\}$
\ \ 
--- a quadratic form is mapped to the set -- ellipsoid of all the vectors
at which this form takes the value 1. Let us explain in what sense ellipsoids
form a smooth manifold. The manifold of ellipsoids is identified with the 
manifold \ $\SO(n) \times R_+^n$ \ by parameterizing it by a pair \
$(g, \spe)$, \ \ $g \in \SO(n)$, $\spe \in R_+^n$. \ Here $g$ is the operator
from $\SO(n)$ which maps the fixed orthonormal basis $(e_1, \ldots, e_n)$ \ 
to the basis of the major axes of the ellipsoid, \ and the tuple \ $\spe$ \
of positive real numbers defines the length of each major axis of the 
ellipsoid. The number of these parameters is \ 
$\dimSO(n) + n = \ (n(n-1)/2) + n = \ n(n+1)/2$.

In this representation, $\cal{M}_+$ \ is the manifold of all the 
{\it ellipsoids of rotation}.

Theorem 1 from Appendix 10 of \cite{Ar1} asserts (in our denotation): 

(1) $\cal{M}_+$ is a finite union of smooth manifolds of codimension 2 
    and higher in the manifold $\Sym_+$,

(2) $\cal{MH}_+$ is a smooth manifold of codimension 2 in $\Sym_+$
\\
(the statement (2) is not in the formulation, but it is in the initial part
of the proof of Theorem 1 in \cite{Ar1}).

Further, by the \ lemma \ in the above Theorem 1 in \cite{Ar1}, \ the set \ \
$\cal{M}_+ \backslash \cal{MH}_+$ \ \ is a finite union of smooth manifolds in
$\Sym_+$ of codimension more than 2.
Due to this, $\cal{MH}_+$ is a set open in $\cal{M}_+$ and dense in 
$\cal{M}_+$.

And by the described above symmetry of scalar shifts, this theorem is also 
true for the surfaces $\cal{M}$ and $\cal{MH}$ in $\Sym(n)$.\\
{\bf The theorem is proved.}

\section{Several definitions and preliminary constructs}
\label{sec-prelim}

\noindent {\bf Definition \ (Basic).}

{\bf 0.} \ $V = \mathbb{R}^n$ \ \ is the Euclidean space.

 $L(v_1, \ldots, v_m)$ \ \ denotes the subspace in $V$ spanned by the vectors
                           $v_1, \ldots, v_m$.

 ${\cal L}(V) = {\cal L}(n)$ \ \ 
                denotes the algebra of the linear operators in $V$\\
                (or --- real square matrices of size \ $n$).

 $E$ \ \ denotes the identity operator in $V$.
\smallskip

{\bf 1.} \ We call an \ irreducible surface \ an irreducible (prime) 
           affine algebraic manifold over any field
           (\cite{Sha}, Book 1, Chapter I, Section 3.1). \ 
  In this paper, all the spaces are over $\mathbb{R}$.\\
  Although in \cite{Sha} all the results are given for an algebraically closed 
  field, our two references to \cite{Sha} are so that the corresponding proofs
  from \cite{Sha} are also valid for the real space.

{\bf 2.} \ $\Dg$ \ denotes the space of diagonal matrices 
   (of dimension \ $n$),

   $\Scal$ \ denotes the line of scalar matrices. 

   Evidently, $\cal{M}$ \ is a conic set containing \ $\Scal$. 
 
   The word ``diagonal'' \ means the main diagonal of a square matrix.

{\bf 3.} \ $A(i,j) = A_{i,j}$ \ are the denotations for the element $a_{i,j}$ 
   of a matrix $A$.\\
   For a diagonal matrix $D$, \ denote \ $D(i) = D(i,i)$.\\
   $\diag(a_1,\ldots,a_n)$ \ denotes the diagonal matrix $D$ having \ 
   $D(1) = a_1, \ldots, \ D(n) = a_n$.
                            
   A \ diagonal-ordered matrix \ is a diagonal matrix $D$ in which \  
   $D(1) \leq \ldots \leq D(n)$.
\smallskip

{\bf 4.} $O(V) = O(n)$ \ \ is (an algebraic, linear) Lie group of orthogonal 
                           operators in $V$.

  $\SO(V) = \SO(n)$ \ \ is the Lie subgroup in \ $O(V)$ \ defined by the 
            condition\\ determinant $= 1$. \ 
            It has dimension \ \ $\dimSO(n) = \ n(n-1)/2$.

  We use certain properties of $\SO(n)$, like \ connectedness, algebraic 
  irreducibility, compactness, smoothness. They are explained in the 
  textbooks on the linear (and/or algebraic) Lie groups. For the Russian 
  language, we give a reference $\cite{VinO}$ for this. But for English, I 
  know only \cite{Bo} (which is of a rather abstract style). 
  Irreducibility of $\SO(n)$ follows from that it is algebraic and connected: 
  \cite{Bo}, head I, Section 1.2, first Proposition.
\smallskip

{\bf 5.} $X \mapsto g \ X \ g^{-1}$  \ is the action of the 
   coordinate (basis) {\it change} by an invertible operator $g$ in a
   matrix (operator) \ $X$. \ \ 
   We shortly call this conjugation action \linebreak ``change by $g$'',
   ``changes by operators from the group $G$'', ``changes from $G$''. 

   And for an orthogonal operator $g$, there holds \ $g^{-1} = g^*$ \ 
   (transposed matrix), \ and \ the expression \ $g \ X \ g^*$ \ 
   represents the basis change in a symmetric bilinear form $X$.

   For \ $g \in O(n)$, \ we denote the action of the basis change as 
   $$
      g^c X = \ g^c(X) = \ g \ X \ g^*,
   $$
   and remember that $g^c$ is an operator in $\Sym(n)$. 
   \ \ \ \ \ \ \ \ \ \ \ \ $g^{-c}$ \ denotes \ $(g^{-1})^c$.

   In the case of dimension 3 \ ($n = 3$), \ \ \ \ 
   1) \ $\RtA(l, \phi)$\\
      denotes the operator of rotation about an axis $l$ by an angle $\phi$, 

   2) \ $l^\phi$ \ denotes the change \ $\RtA(l,\phi)^c$, 

   3) \ 1-orbit of a form $M$ by an axis $l$ \ means the orbit of $M$ 
      under the changes by the rotation operators about $l$.
\smallskip

{\bf 6.} Denote \ $\cal{MD}$ \ the set of diagonal matrices in $\cal{M}$, \ \ 
   $\cal{MDH} = \cal{MD} \cap \cal{MH}$.
   \\
   For $1 \leq i < j \leq n$ \ denote \ $\Pi_{i,j}$ \ the plane in $\Dg$ 
   defined by the equation \ $D(i) = D(j)$. 

   For example, for \ $n = 3, \ \ \Pi = \Pi_{1,2}$ \ is all the matrices
   of the kind \ $\diag(\lambda,\lambda,\mu)$, \ while \ $\Pi_{1,3}$ and \ 
   $\Pi_{2,3}$ are obtained from $\Pi$ by permutational changes from $\SO(3)$. 

   Denote \ $\Pi_{i,j}^h = \ \Pi_{i,j} \cap \cal{MH}$.
\smallskip

{\bf 7.} In the sequel the word ``orbit'' means: \ the orbit \ $\Orbit(A)$ \ 
     of a form \ $A \in \Sym(n)$ \ under the action of changes from $\SO(n)$.

{\bf 8.} Fix a \ scalar product (Euclidean metric) \ $\scp$ \ on the space \ \ 
  $V = \mathbb{R}^n$, \ \ an orthonormal by $\scp$ basis \ \ 
  $\Bas = \{e_1, \ldots, e_n\}$ \ \ in \ $V$, \ \ 
  and let the matrices in ${\cal L}(n)$ represent the linear operators in $V$ 
  written in \ $\Bas$. \ \ 
  Also let the matrices in \ $\Sym(n)$ \ represent the symmetric bilinear 
  forms on \ $V$ \ written in \ $\Bas$.

  We also speak of operators from $O(n)$ given by permutations on the set 
  $\Bas$, or, for example (in the case of $n = 3$), being rotations about
  an axis $e_i$.

{\bf 9.} \ $\trace(X)$ \ denotes the sum of the diagonal elements of a matrix 
           $X$.    

{\bf 10.} \ The \ stabilizer \ $\St(A)$ \ of a form $A \in \Sym(n)$ \ \  
     is the subgroup of operators $g$ in a group ($O(n)$ or $\SO(n)$) which 
     change preserves $A$ (that is \ $g A = A g, \ g^c A = A$).

{\bf 11.} \ The \ eigenspace \ of a given eigenvalue $\lambda$ of a linear 
   operator $A$ in a space $V$ \ is the subspace of all the vectors $v$ 
   such that \ $A v = \lambda v$.

{\bf 12.} \ For a commutative ring $R$ and a set $S \subset R^n$, \
  a map \ $P : S \rightarrow R^m$\\  is called \ {\it polynomial}, \ if
  each coordinate of its image is represented as a polynomial in the 
  coordinates $\{x_1, \ldots, x_n\}$ of the support domain.

  For example, the map \ $P : \ \Dg \times \SO(n) \rightarrow \Sym(n)$ \ \  
  given by the formula \ $P(D, g) = g^c \ D$ \ is {\it polynomial}, because 
  each element of the result matrix is expressed as a polynomial in 
  the elements of the matrices $D$ and $g$.
\smallskip

\noindent {\bf Example:} \ 
for $n = 3$, \ the diagonal matrices in $\cal{MH}$ are exactly the matrices 
of the kind \ $\diag(\lambda,\lambda,\mu)$ \ with $\lambda \neq \mu$, \  
and also the two families obtained from this one by permutations on the 
diagonal.
\smallskip

Our discourse bases on \ known facts of the linear Lie group theory \ 
and the following classical 
\smallskip

\noindent {\bf Theorem DO.} \ \ 
(1) For any (real) quadratic form $M$ there exists a change \ 
    $g \in O(n)$ \ which brings $M$ to a (real) diagonal matrix having on 
    the diagonal all the eigenvalues of $M$.

Orthogonality in this theorem is considered with respect to the scalar
product \ $\scp$ \ introduced in the above definition (Basic).

(2) Each eigenvalue $\lambda$ for $M$ defines the non-zero eigenspace of
    $\lambda$ for $M$.

(3) Eigenspaces belonging to different eigenvalues of $M$ are orthogonal to 
    each other.

(4) (our addition) \ For bringing $M$ to such a diagonal matrix there are 
    sufficient the changes from $\SO(n)$. 
\smallskip

This theorem is proved, for example, in \cite{Wa}, Head XII, \S 90. The
statement (1) is by the theorem in the end of \S 90 about a pair of quadratic 
forms, and we put one of these forms to be \ $\scp$. \ 
The proof for the points (2) and (3) is contained in the discourse in the 
whole \S 90.

But it remains to add the proof for the point (4).
Let $h$ be the operator in $O(n)$ defined by transposition of the vectors 
$e_1$ and $e_2$, \ and let \ $g$ bring $M$ to a diagonal matrix. 
If \ $\det(g) = 1$, \ then the statement is true. Otherwise, $\det(g) = -1$, \ 
and the composition \ $h g$ \ belongs to $\SO(n)$ and still brings $M$ 
to a diagonal form.
\smallskip

\noindent {\bf Remark.} \ \ 
As \ $\charPol(A)$ \ does not change with any invertible coordinate change in 
$A$, and due to Theorem DO, \  
the spectrum width for a matrix in $\Sym(n)$ is the spectrum width of its 
diagonal form.
\smallskip

\noindent {\bf Lemma Stab.} \ \ 
Let \ $M \in \cal{MH}$. 

1) The set $V_1$ of the vectors belonging to the multiple eigenvalue
   for $M$ is a subspace of dimension 2.
 
2) The connected component $G_e$ of unity in the stabilizer $\St(M)$ 
   in $\SO(n)$ is a single-parameter subgroup of the operators which 
   restriction on $V_1$ are rotations of this 2-dimensional space about 0, 
   and which restriction to the orthogonal complement to $V_1$ is identity.

Appendix Stab \ provides a simple proof for this.
\smallskip

\noindent {\bf Definition s-metric.} 
\ \ 
Define a scalar product \ $\langle_s, \ \rangle$ \ on the space 
${\cal L}(n)$: \ 

\centerline{
    $\langle_s X, Y\rangle = \sum_{1 \leq i,j \leq n} x_{i,j} y_{i,j}$,
}\noindent
and call it \ s-metric.\\
It is an Euclidean metric: a positive definite bilinear function.\\
Denote \ \ $\sQuad(X)$ \ \ its corresponding quadratic form --- 
sum of squares of the elements of a matrix $X$.

Note that s-metric defines the same topology in ${\cal L}(n)$ as the usual 
linear operator norm.

\noindent \underbar{Agreement}\\
In this paper, the expressions \ \ 
``distance'', ``angle'', ``orthogonality'', ``isometry'',\linebreak
``orthogonal projection'', ``straight cylinder'', ``circumference'', 
``sphere'', ``circumference center'', ``radius'' \ \  
applied to points and subsets in ${\cal L}(n)$ \ --- are understood in the 
sense of s-metric\\ 
(do not confuse this with the metric \ $\scp$ \ of the space \ $V$ \ on 
which the operators from ${\cal L}(n), \Sym(n)$ act!).
\smallskip

\noindent {\bf Lemma S} \ (probably, known).
  
  For any \ $X \in {\cal L}(n)$ \ and any orthogonal operator $g$, 
  it holds the law

  \verb#             # $\sQuad(X) = \sQuad (g X) = \ \sQuad (g^c X)$.

In particular, any change from $O(n)$ is an isometry of ${\cal L}(n)$ 
in the s-metric.

\noindent {\bf Proof.} \
By the definition of an orthogonal operator $g$, such an operator preserves 
the scalar square for each column--vector $C_j$ in $X$: \  
$\sQuad(g C_j) = \sQuad(C_j)$. \ As the $\sQuad(X)$ is sum of the 
scalar squares of the columns in $X$, it holds \ $\sQuad(g X) = \sQuad(X)$. \
Similarly, the right-hand side multiplication by $g^*$ preserves the scalar 
square of each row, and hence the s-square of $X$.  
Therefore \ $\sQuad (g^c X) = \ \sQuad(X)$. \ \ 
{\bf The Lemma is proved.}
\smallskip

\noindent {\bf Definition.} 
 
 1) A {\it shift} at a matrix $A$ \ is the map of adding of $A$ in the 
    affine space ${\cal L}(n)$.

 2) A {\it scalar shift} \ is a shift at a scalar matrix.

\noindent
We use that a shift at a matrix is an isometry of the affine space of 
${\cal L}(n)$. 
This is because a shift at a vector in an affine space preserves the vector 
defined by any pair of points.

\section{Preliminary notions about orbit}
\label{sec-OrbitPrelim}

\noindent {\bf Lemma Orbit1.} 
 \ \ 
 (1) Any change from $\SO(n)$ maps $\cal{M}$ on itself.

 (2) Scalar matrices are the only fixed points of the action of $\SO(n)$ 
     by changes on $\Sym(n)$.

 (3) $\cal{M}$ is union of mutually non-intersecting orbits.

 (4) Each matrix in $\cal{MD}$ is brought by some change from $\SO(n)$ to an 
     ordered diagonal. \ 
     In particular, each permutation at the diagonal of a matrix in 
     $\cal{MD}$ is represented by some change from $\SO(n)$.

 (5) Two forms in $\cal{M}$ belong to the same orbit if and only if
     they have the same spectrum.

 (6) For each matrix in $\cal{M}$, its orbit contains exactly one 
     diagonal-ordered matrix.

 (7) 7.1. For each real number \ $t$ \ and a form \ \ $A \in \Sym(n)$, \ \ 
          it holds\\ $\Orbit(t \cdot A) = \ t \cdot \Orbit(A)$. 

     7.2. For \ $n = 3$, \ the orbits of any two forms in $\cal{MH}$ are 
          obtained from each other by some scalar shift and a homothety by 
          some non-zero factor\\
          (and we use that this transformation preserves the angle between
           any two straight lines in $\Sym(n)$).

 (8) 8.1. Discriminant of a matrix in ${\cal L}(n)$ does not change under a 
          scalar shift. Any scalar shift maps the surface $\cal{M}$ on itself.

     8.2. For each real number \ $s$, \ intersection \ $\cal{M}$$_s$ \ of \
          $\cal{M}$ \ with the hyper-plane\\ $\trace(X) = s$ \ \ 
          is a surface mapped on itself by any change from $\SO(n)$. 

     $\cal{M}$ \ is a straight cylinder over \ $\cal{M}$$_0$ \ having the 
     line \ $\Scal$ \ as element.

\noindent {\bf Remark.}\\ 
In the statements (4), (5), (6) it is essential that a matrix has a multiple 
eigenvalue.
\smallskip

\noindent {\bf Proof}. 
The statements (1), (2) and (3) are known and evident.

\noindent \underbar{Proving (4)} \  
$D$ has at least two equal elements on the main diagonal, let them correspond
to the basis vectors $u$ and $v$ respectively. Then the operator $g_2$ which 
transposes $u$ and $v$ and is identical on the rest of the basis, has 
determinant -1. \ Some permutation $\sigma$ on the basis brings the main 
diagonal to an {\it ordered diagonal} 
(definition (Basic) in Section \ref{sec-prelim}) 
and is represented by some operator $g \in O(n)$. If $\sigma$ is even, 
then $g$ satisfies the 
requirement of the point (4). Otherwise, the composition \ $g \ g_2$ \  
belongs to $\SO(n)$, and its change brings $D$ to a diagonal-ordered form. \ \
{\it The point (4) is proved}.

\noindent \underbar{Proving (5) and (6)}. \ 
The change action preserves the coefficients of the characteristic 
polynomial, so, it preserves the spectrum. \ 
Let us prove the reverse part of the statement (5). 
Let $A$ and $B$ have the same spectrum. By Theorem DO, $A$ and $B$ are 
diagonalized by some changes from $\SO(n)$. By the statement (4), they are 
brought further to diagonal-ordered matrices. As these diagonal-ordered 
matrices represent the same multi-set, they are equal. \ 
{\it And this also proves the statement} (6). 

\noindent \underbar{Proving (7)} 
\ \
The statement (7.1) is evident. Let us prove (7.2).

By Theorem DO, the orbit of any matrix in $\cal{MH}$ coincides with the 
orbit of some diagonal matrix in $\cal{MH}$, and by the statement (5),
it is the orbit of some matrix \ $D_1$ \ of the kind \  
$\diag(\lambda, \lambda, \mu)$ \ with $\lambda \neq \mu$. \ 
Therefore, it is sufficient to prove (7.2) for any diagonal matrices \ 
$D_1, D_2 \in \cal{MH}$, \ where $D_1$ is given above, \  
$D_2 = \diag(\lambda_2, \lambda_2, \mu_2), \ \ (\lambda_2 \neq \mu_2)$.

The matrix \ $D = \diag(0, 0, 1)$ \ is transformed to $D_1$ by a shift at 
some matrix of the kind \ $C = \diag(a,a,a)$ \ \ and by a homothety by 
some non-zero factor \ $b$: \ \ $(D+C) \cdot b \ = D_1$. \ \
For such representation, it is sufficient to put \ \ 
$a = \lambda/(\mu - \lambda)$, \ \ $b = \mu - \lambda$. \ \ 
Each change $g^c$ is an isomorphism on the algebra of square matrices. 
In particular, for each square matrix $X$ there hold the equations
$$
  g^c (X+C) = (g^c X) + (g^c C) = (g^c X) + C, \ \ \ 
  g^c (b \cdot X) = b \cdot g^c X.
$$
\noindent
Therefore, with adding a scalar matrix $C$ to any symmetric matrix $X$,
the orbit shifts at $C$, \ and with multiplying $X$ by a coefficient $b$, 
the orbit is homothetically multiplied by $b$.
Therefore, $\Orbit(D_1)$ is obtained from $\Orbit(D)$ by an appropriate
scalar shift and a homothety required in the point (7.2).

In the same way, by a scalar shift and an invertible homothety, $D_2$ is
converted to $D$. The composition of these shifts and homotheties gives
the transition \ $D_2 \mapsto D \mapsto D_1$. \ 
And such a transformation equals to composition of a single scalar shift 
and a single invertible homothety --- due to the equations

  $(((X + a \cdot E) \cdot b) + a_2 \cdot E) \cdot b_2 = \                     
    b b_2 ((X + a \cdot E) + (a_2/b) \cdot E) =$

  $b b_2 ((X + (a + a_2/b) \cdot E)$.
\\ 
{\it The statement (7) is proved}.

\noindent \underbar{Proving (8)} \ \ 
{\bf 8.1:} adding of a scalar matrix \ $c \cdot E$ \ shifts the spectrum of 
  any matrix $A$ in ${\cal L}(n)$ at the number $c$.  This follows from the 
  definition of a vector belonging to an eigenvalue $\lambda$: \ \ 
$A v = \lambda v \ \ ==> \ \ (A + c E) v = \lambda v + c v = (\lambda + c) v$
\ \
(note that for a non-symmetric matrix the eigenvalues and their vectors 
may be not real).

This shift preserves discriminant, because (by the discriminant definition via
the eigenvalues) it depends only on the pairwise differences of eigenvalues.\\
{\it The statement (8.1) is proved}.

{\bf 8.2:} \ 
recall that the surface $\cal{M}$ resides in $\Sym(n)$ and contains
the line $\Scal$ of scalar matrices. \ 
For each real number $s$, denote \ $\Pi(s)$ \ the hyper--plane in $\Sym(n)$ 
defined by the equation \ $\trace(X) = s$, \ \  
$\cal{M}$$_s$ $= \cal{M} \cap$ $\Pi(s)$. \ \ So, \ 
$\cal{M} = \cup \{\cal{M}$$_s | \ s \in \mathbb{R}\}$.

Any change from $\SO(n)$ preserves the matrix trace, as well as the 
discriminant value. Therefore, such a change maps each plane $\Pi(s)$ on 
itself and maps each restriction $\cal{M}$$_s$ on itself. 
So, for the point 8.2, it only remains to prove the statement about the 
straight cylinder.

Parameterize the line $\Scal$ by the value $s$ of a matrix trace: \ 
$s \mapsto C(s)$. \ Use also that $\Scal = L(E)$.
It follows from the s-metric definition that for all 
$X \ \langle_s X, E\rangle = \trace(X)$. \ 
For each $X$ from $\Pi(s)$ \ $X$ and $C(s)$ have the same trace, hence 
the vector \ $\overline{C(s) X}$ \ is a matrix of zero trace. 
Hence \ $\langle_s \overline{C(s) X}, \ E\rangle = 0$ \ for all 
$X \in \Pi(s)$, \ and the hyper-plane $\Pi(s)$ is orthogonal to $\Scal$. \ 
As the surface $\cal{M}$$_s$ is in $\Pi(s)$, the line $\Scal$ orthogonally
intersects $\cal{M}$$_s$ (in the point $C(s)$).  

Each fiber $\cal{M}$$_s$ is obtained from $\cal{M}$$_0$ by the scalar shift 
at $C(s)$, in particular, all these fibers are isometrical to $\cal{M}$$_0$.
Summing up the discourse, conclude that\linebreak 
1) $\cal{M}$ contains the line $\Scal$, \ \ 
2) $\cal{M}$ is union of the fibers, \ \ 
3) each fiber is obtained from $\cal{M}$$_0$ by the scalar shift at $C(s)$, 
   intersects with $\Scal$ in $C(s)$, and is orthogonal to $\Scal$.  

This means that $\cal{M}$ is the straight cylinder over $\cal{M}$$_0$
having the line $\Scal$ as element.\\
{\bf The lemma is proved}. 
\medskip

\noindent {\bf Lemma A0} (probably, known). \ 

Commutator of a diagonal and an anti-symmetric square matrices 
(of the same size) is a symmetric matrix having zero diagonal.

\noindent {\bf Lemma A1} (known). \
The counter-image $\cal{M}^\prime$ of an algebraic set $\cal{M}$ for a 
polynomial map \ $F : \ R^n \rightarrow R^m$ \ is an algebraic set.

\noindent {\bf Lemma A2} (probably, known). \ 
For a surjective polynomial map from an algebraic set $\cal{M}$ onto an 
algebraic set $\cal{M}^\prime$, \ if $\cal{M}$ is irreducible, then 
$\cal{M}^\prime$ is irreducible.
\smallskip

\noindent Simple proofs for lemmata A0, A1, A2 \ are given in Appendix.
\smallskip

\noindent {\bf Remark.} \ 
The lemmata A1 and A2 \ are related to the topology of the algebraically 
closed sets in an affine space. The corresponding notions are given, for 
example, in \cite{Sha}, Head I, \S 2, \S 3.

\newpage
\section{Theorem of irreducibility}
\label{sec-prime}

\noindent {\bf Theorem ``Irreducibility''.}
 
The discriminant surface $\cal{M}$ is algebraically irreducible

(that is it is not a union of any two algebraic sets different from 
$\cal{M}$).  
\medskip

\noindent {\bf Proof.}\\ 
Let \ \ $\Pi = \mathbb{R}^{n-1}$ \ \ be the plane parameterized by the tuple \ 
$(\lambda,\mu_1,\ldots,\mu_{n-2})$ \ of real numbers. 
Consider the map \ \ 
 $\DO : \Pi \times \SO(n) \rightarrow \cal{M}$, 

\verb#                      # $\DO ((\lambda,\mu_1,\ldots,\mu_{n-2}), \ g) = \ 
                            g^c \diag(\lambda,\lambda,\mu_1,\ldots,\mu_{n-2})$.
\smallskip

\noindent
This is a polynomial map of algebraic surfaces (where the support domain has 
one dimension more than the image).

Let us prove that this map is surjective. Let $A \in \cal{M}$. 
By Theorem DO, the orbit of $A$ contains some diagonal matrix $D$. 
By the statement 5 of Lemma Orbit1, there exists a multiple eigenvalue in the 
spectrum of $D$ (and $A$), denote it \ $\lambda$. \ 
By the statement 4 of Lemma Orbit1, each permutation on the diagonal of $D$ 
is represented by some change from $\SO(n)$.
And some permutation produces a matrix of the kind \ 
$\diag(\lambda,\lambda,\mu_1,\ldots,\mu_{n-2})$. \ Hence, \ 
$A = \DO((\lambda,\mu_1,\ldots,\mu_{n-2}), \ g)$ \ for some $g \in \SO(n)$. 
\ \ 
{\it The surjectiveness is proved}.

As it is mentioned in the beginning of Section \ref{sec-prelim}, \ 
$\SO(n)$ is an irreducible algebraic surface.
The plane $\Pi$ \ is also irreducible. 

{\it Direct product of any two irreducible algebraic sets is irreducible.} \ \
To prove this can be an exercise. Also the proof is given, for example, in 
\cite {Sha}, Head I, \S 3.1. \ 
By this theorem, the algebraic surface \ $\Pi \times \SO(n)$ \ is irreducible. 
Hence, by Lemma A2, \ $\cal{M}$ is irreducible.\\ 
{\bf The theorem is proved}.

\newpage
\section{Structure of the discriminant surface in $\Sym(n)$}
\label{sec-DSForm}

Below, $\Pi$ denotes the plane $\Pi_{1,2}$ of dimension \ $n-1$ \ in the 
surface $\cal{MD}$ \ --- see Definition (Basic). \ Below we identify
it to \ $\mathbb{R}^{n-1}$ \ by the parameterization  
\ \ \ $(\lambda, \mu_1, \ldots, \mu_{n-2}) \ \mapsto \    
       \diag(\lambda, \lambda, \mu_1, \ldots, \mu_{n-2})$.

Also keep in mind that \ in the plane $\Pi$ the matrices from $\cal{MH}$ form
an open set.

Recall also that the word ``orbit'' means: \ the orbit of a quadratic form
under the action of changes from $\SO(n)$.
\smallskip

\noindent {\bf Theorem DS} 

{\bf (1)} $\cal{M}$ is an irreducible algebraic variety.

{\bf (2)} $\cal{M}$ \ is union of the restrictions \ $\cal{M}$$_s$, \ 
          each defined by the equation \ $\trace(X) = s$.

   $\cal{M}$ is the straight cylinder over $\cal{M}$$_0$ having the line 
             $\Scal$ as element.

{\bf (3)} Orbits for the matrices in $\cal{M}$ are classified by spectrum,
          or --- by the unique diagonal--ordered matrix in the orbit.

{\bf (4)} $\cal{MH}$ has the following structure.

 {\bf 4.1.} Locally, $\cal{MH}$ \ is a smooth fibration over the plane $\Pi$ 
  with the fiber being an orbit of a diagonal matrix from $\cal{MH}$.

  The orbit of each matrix in $\cal{MH}$ orthogonally intersects the plane 
  $\Dg$ of diagonal matrices.

 {\bf 4.2.} Orbit of each matrix in $\cal{MH}$ is \ a smooth, compact, 
  algebraic, connected surface of dimension \ 
  $\dimSO(n) - 1 = (n(n-1)/2) - 1$. \ It resides on a sphere in the
  hyper-plane orthogonal to the line of scalar matrices. 

{\bf (5)} The map \ $\DO : \Pi \times \SO(n) \rightarrow \cal{M}$,
 
    \hspace{25mm}   $\DO((\lambda,\mu_1,\ldots,\mu_{n-2}), \ g) = \ 
                     g^c \diag(\lambda,\lambda,\mu_1,\ldots,\mu_{n-2})$
\\
 is a surjective polynomial map of algebraic surfaces (where the
 support domain has one dimension more than the image).
\medskip

\noindent {\bf Remark.}
 
{\bf 1.} \ A local parameterization of $\cal{MH}$ is the same as a 
  parameterization of $\cal{M}$ in some neighbourhood of a point in 
  $\cal{MH}$. \ Because by Theorem Ar (Section \ref{sec-Ar}), $\cal{MH}$ 
  is open in $\cal{M}$.

{\bf 2.} The statement (4.1) \ does not differ much from \ Theorem 1 \ 
from \cite{Ar1} \ about smoothness of $\cal{MH}$ in the area of ellipsoids 
(Section \ref{sec-Ar} of this paper, Theorem Ar, point (2)).
Therefore, we skip the proof for this point, instead we provide the following
informal description for the local parameterization near a point in 
$\cal{MH}$. \ Due to the action of $\SO(n)$, the local structure of $\cal{MH}$
near each point is isometrical to a neighbourhood of a diagonal matrix 
$D_0 \in \cal{MH}$. In some neighbourhood $U$ of $D_0$, the surface 
$\cal{MH}$ is a smooth fibering to orbits, where each fiber-orbit comes out 
from \ $D \in \Pi \cap U$. \ 
The orbit of $D$ is locally diffeomorphic to the homogeneous space \ 
$\SO(n)/\St(D)$ \ of the congruence classes by the stabilizer. 
And by Lemma Stab, the connected component $St_1$ of unity in the 
stabilizer $\St(D)$ is a single-parametric subgroup of rotations of the 
2-dimensional eigenspace $V_1$ of the (unique) multiple eigenvalue in $D$ \ 
(taking in account that the operators from $St_1$ are identical on the 
orthogonal complement to $V_1$).

The used here notions about a \ stabilizer, orbit, and a smooth structure 
on the space of congruence classes, \ we find in the book \cite{VinO} 
(in Russian), head I, Theorem 1, and the following exercise; \ 
Theorems 3, 4, section 9 about homogeneous spaces.\\ 
For English, the corresponding theory is, hopefully, given in \cite{Bo}.

In terms of the manifold of ellipsoids (\cite{Ar1}, Appendix 10), this 
statement (4.1) about a smooth fibration is expressed as follows.
In the parameterization of a neighbourhood of a matrix $D \in \mathbb{MDH}$ 
by the pair \ $(g, \spe)$ (see Section \ref{sec-Ar}), the part \ $\spe$ \ 
corresponds to the plane $\Pi$ in the space of diagonal matrices. And under
infinitesimal rotations $g$ (close to $E$) of the basis of the major 
ellipsoid axes, each axis turns in a direction orthogonal to this axis.

\section{Addition to the theorem for the case of $n = 3$}
\label{sec-DS3Form}

Recall that \ $\cal{M}$$_s$ = $\{A \in \cal{M} |$ \ $\trace(A) = s \}$. 
\smallskip

\noindent {\bf Theorem DS3.}

{\bf (1)} Singular points on $\cal{M}$ are the scalar matrices and only them.

{\bf (2)} $\cal{M}$ is the straight cylinder over $\cal{M}$$_0$ 
 (Theorem DS (2)). Therefore, it is sufficient to describe the surface 
 $\cal{M}$$_0$. This 3-dimensional surface resides in a 5-dimensional
 space and has the following structure.

{\bf 2.1.} $\cal{M}$$_0$ \ is the cone with zero as vertex over the orbit 
  of the matrix \ $\diag(1,1,-2)$. This 2-dimensional orbit resides on a 
  4-dimensional sphere having center in zero.

{\bf 2.2.}
  The orbit of each non-scalar matrix in $\cal{M}$ is a 2-dimensional
  algebraic surface residing on a 4-dimensional sphere. 
  It is diffeomorphic to the projective plane $\mathbb{R}P(2)$.

{\bf (3)} The orbit of each matrix in $\cal{M}$ is obtained from the orbit of
          the matrix\linebreak $D = \diag(1,1,-2)$ \ by some scalar shift and 
          some homothety.
 
{\bf (4)} The set of diagonal matrices in \ $\cal{M}$$_0$ \ consists of the 
  three straight lines which intersect in zero and can be parameterized
  respectively as\\ $\diag(\lambda, \lambda, -2 \lambda)$, \ \ 
                    $\diag(\lambda, -2 \lambda, \lambda)$, \ \ 
                    $\diag(-2 \lambda, \lambda, \lambda)$.
\smallskip

\noindent {\bf Remark.} 
A natural smooth parameterization of $\cal{M}$$_0$ in some neighbourhood of \ 
each form $0 \neq A$ in $\cal{M}$$_0$ \ is the fibering over the line 
$\overline{0 A}$ with the fiber being the orbit of a matrix from this line, 
and with this orbit intersecting orthogonally this line. 
\smallskip

All the remaining discourse serves the proof of the theorems DS, DS3.

\newpage
\subsection{Investigating an orbit by orthogonal changes} 
\label{sec-DSForm:orbit}

\noindent {\bf Lemma (AlgSO)} \ 
{\bf (1)} The tangent space to the surface \ $\SO(n)$ \ in the point \ $E$

     (the Lie algebra $\so(n)$) consists of all the antisymmetric matrices.

{\bf (2)} The tangent space to the orbit of a diagonal matrix $D$ in the 
  point $D$ \ consists of the commutators of $D$ with all the antisymmetric
  matrices.

This lemma \ a) is known, \ b) is an exercise for the reader, \ 
c) is simply derived in Appendix (AlgSO).

\noindent {\bf Lemma (Orbit2).} 

{\bf (1)} The orbit of any non-scalar matrix in $\Sym(n)$ \ 
          intersects orthogonally the plane of {\it diagonal} matrices.

{\bf (2)} Orbit of any matrix in $\Sym(n)$ resides in the hyper-plane  
          orthogonal to the line $\Scal$.

{\bf (3)} Orbit of any matrix in $\cal{MH}$ is a \ smooth, compact, 
          algebraic, connected surface of dimension \ \ $\dimSO(n) - 1$.

\noindent {\bf Proof.} 

\noindent \underbar{Proving (1)} \ 
Due to the symmetry by the group action, and by Theorem DO, it is 
sufficient to prove this statement for the orbit of any diagonal matrix $D$.
By Lemma AlgSO, the tangent space $T$ to the orbit of $D$ in the point 
$D$ consists of commutators of $D$ with all antisymmetric matrices.
Hence, by Lemma A0, each matrix $M$ in $T$ has zero diagonal. Then, it folows 
from the definition of the s-product (Section \ref{sec-prelim}) that \ 
$\langle_s D_1, \ M\rangle = 0$ \ for each diagonal matrix $D_1$. \ \
{\it The statement (1) is proved.}

\noindent \underbar{Proving (2)} \ \ 
Let $A \in \Sym(n)$. \ If $A$ is scalar, then its orbit is $\{A\}$, and the
statement (2) is satisfied. It remains the case of a non-scalar $A$. \ 
Denote $\Pi^*$ the hyper-plane orthogonal to $\Scal$ and containing $A$. 
Let $A_0$ be the base of this projection of $A$, that is the intersection
point of $\Pi^*$ with $\Scal$.
By Lemma S, each orthogonal change is an isometry in $\Sym(n)$. And each point 
in $\Scal$ is a fixed point with respect to these changes. Hence, as the 
vector $\overline{A_0 A}$ is orthogonal to $\Scal$, the vector 
$\overline{A_0 \ (g^c A)}$ is also orthogonal to $\Scal$ for each 
$g \in SO(n)$. \ \ 
{\it This proves the statement (2).}

\noindent \underbar{Proving (3)} \ 
First, let us present a finite polynomial system defining the orbit of a 
matrix $A$ from $\cal{M}$. Let \ $c_i(X)$ \ be the coefficient in the 
monomial of $\lambda^i$ in the polynomial $\charPol_X(\lambda)$ for an 
arbitrary symmteric matrix $X$ \ ($c_i$ is a polynomial in the variables 
$x_{i,j}$ --- elements of $X$). The value $c_i(A)$ for each concrete symmetric
matrix $A$ (consisting of real numbers) is expressed by the Viete formula
via \ $\spec(A)$:

\centerline{$c_i(A) = \sgn(i) \cdot \elSym_i(\ev(A)), \ \ 0 \leq i \leq n-1$.
}
\noindent Here \ $\ev(A)$ \ is the number sequence made of \ $\spec(A)$, \ 
with repetition of multiple values,

$\elSym_i$ \ is the elementary symmetric polynomial of degree \ $n-i$ \ in 
  the variables \ $y_1, \ldots, y_n$ \ (in this formula the real eigenvalues
  are substituted for the variables $y_j$ in $\elSym_i$),

$\sgn(i)$ \ is the appropriate sign: plus or minus.

Therefore, $\Orbit(A)$ is defined by the following equations for unknown 
matrix $X \in \Sym(n)$:

\centerline{$c_i(X) = \sgn(i) \cdot \elSym_i(\ev(A)), \ \ 0 \leq i \leq n-1
                                               \ \ \ \ \ \ \ \ \ \ \ \ (1)$.
}\noindent 
Here each lefthand side is the fixed polynomial $c_i$ 
(in the variables $X$) defined above, and the righthand side is a real 
number --- the value of the Viete expression at the spectrum of $A$.

Let us prove that the system (1) defines the set $\Orbit(A)$. \ 
If $X \in \Orbit(A)$, \ then \ $\spec(X) = \spec(A)$, \ this defines the 
characteristic polynomial for $X$, and $X$ satisfies (1). \ 
Reversely, if for $X \in \Sym(n)$ the equations (1) on $X$ are satisfied, 
then \ $\charPol(X) = \charPol(A)$, \ and so, \ $\spec(X) = \spec(A)$. \ 
By the statement 5 of Lemma Orbit1 of the Section \ref{sec-prelim}, 
the orbit of $A$ consists exactly of the forms in $\cal{M}$ which spectrum 
is \ $\spec(A)$. Hence, $X \in \Orbit(A)$.
\\
{\it The statement about algebraic equations for the orbit is proved}.
\smallskip

\noindent
\underbar{Now, prove the remaining part of the statement (3) of the lemma}

\subsubsection{Further discourse on the orbit}
\label{sec-DSForm:orbit:further}

The map \ \ \ \ \ 
        $\DO_A : \SO(n) \rightarrow \Orbit(A)$, \ \ \ $\DO_A(g) = g^c A$
\\
is a polynomial map of algebraic surfaces. By the definition of orbit, 
$\DO_A$ is surjective. As $\SO(n)$ is compact and connected, and $\DO_A$
is continuous, \ $\Orbit(A)$ is compact and connected. 

Now, irreducibility of $\Orbit(A)$ follows from 

a) that $\SO(n)$ is irreducible (\cite{Bo}, head I, Section 1.2, first 
   Proposition), 

b) that $\DO_A$ is a surjective map onto $\Orbit(A)$, \ \ \
c) Lemma A2.

So: $\Orbit(A)$ is an irreducible, connected, compact, algebraic surface
contained in some sphere in the hyper-plane $\Pi^*$. The minimal dimension
of embracing plane for $\Orbit(A)$ depends on the multiplicities in \ 
$\spec(A)$.

Let us describe a smooth parameterization of the orbit for the form 
$A \in \cal{MH}$. \ Let \ $\Bas = \{e_1, \ldots, e_n\}$ \ be an eigenvector 
basis for $A \in \cal{MH}$, and let $e_1$ and $e_2$ belong the (unique) 
multiple eigenvalue for $A$. \ By Lemma Stab, the connected component of
unity in the stabilizer $\St(A)$ is the group of rotations about 0 of the 
plane $V_1 = L(e_1, e_2)$ \ 
(taking in account that each operator from this component is identical on 
the orthogonal complement to $V_1$). 

$\St(A)$ is a Lie subgroup of dimension 1, and a homogeneous space \ 
$\SO(n)/\St(A)$ \ of the left congruence classes by $\St(A)$ is a smooth
manifold of dimension \ $\dimSO(n) - 1$. \ 
Respectively, $\Orbit(A)$ is a smooth and closed sub-manifold is $\Sym(n)$
of dimension \ $\dimSO(n) - 1$ \ locally diffeomorphic to \ 
$\SO(n)/\St(A)$. \ 
The theory for these three latter statements about the stabilizer and 
the orbit is referenced in the Remark 2 after Theorem DS.

Due to the action of $\SO(n)$, small neighbourhoods of different points
in $\Orbit(A)$ are isometrical. More precisely, for each points $A$ and $B$ 
in this orbit there exists a neighbourhood of $A$ in the orbit which is mapped
onto the corresponding neighbourhood of $B$ in the orbit by the change of 
some operator from $\SO(n)$. This provides a homogeneous smooth atlas 
parameterizing the orbit.\\
{\bf Lemma Orbit2 is proved.} 
\smallskip

{\bf So, Theorem DS is proved. \ The proof parts are:}

(1): irreducibility (Theorem ``Irreducibility''), 

(2): the structure of $\cal{M}$ of a straight cylinder 
     (Lemma Orbit1, point 8),

(3) --- Lemma Orbit1, point 5,

(4.1): \ the orbit of any matrix in $\cal{MH}$ orthogonally intersects
         the space of diagonal matrices (Lemma Orbit2, point 1),

(4.2) --- Lemma Orbit2,

(5) (surjectiveness of $\DO$) --- in the proof of Theorem ``Irreducibility''.

\subsection{Orbit structure details in the case of $n = 3$}
\label{sec-DSForm:orbit3}

\noindent {\bf Lemma (ST).} \ \ 
In the case of \ $n = 3$, \ the following holds for the orbit of any 
non-scalar matrix $A$ in $\cal{M}$.

{\bf (1)} It is a bi-dimensional, algebraic surface residing on a 
          4-dimensional sphere and diffeomorphic to the projective plane 
          $\mathbb{R}P(2)$.

{\bf (2)} It contains the union $O_{cir}$ of a smooth single-parameter family 
          of circumferences of the same radius having a common point $A$\\
          (we guess that $\Orbit(A) = O_{cir}$, but do not prove, so far).

{\bf (3)} It is not contained in any four-dimensional plane.

\noindent {\bf Proof.} 

The statement of that $\Orbit(A)$ is algebraic is proved in Lemma Orbit2.

By the point 2 of Lemma Orbit2, \ $\Orbit(A)$ resides in the 5-dimensional 
plane --- the orthogonal complement to the line $\Scal$. 
The action of $\SO(V)$ preserves the s-metric (Lemma S).
This adds the equation of a sphere for an arbitrary matrix in $\Orbit(A)$,
and it turns out that this orbit resides on a four-dimensional sphere
(by the s-metric).

Let us prove the statement about a diffeomorphism from $\mathbb{R}P(2)$ 
onto the orbit. 
Due to the action of $\SO(3)$, and by the statement (7) of Lemma Orbit1, it 
is sufficient to prove this for the diagonal matrix \ $D = \diag(1,1,-2)$ \ 
(of zero trace).

Represent $\mathbb{R}P(2)$ as the smooth variety of axes in the 3-dimensional 
space $V$, and define the map 

\centerline{$s : \ \mathbb{R}P(2) \rightarrow \ \Orbit(D), \ \ \ 
            s(l) = \Rt(e_3, l)^c D$.
}\noindent
Here \
1) the axis $e_3$ belongs the eigenvalue \ -2 \ of $D$,

2) $l$ \ is an arbitrary axis in $V$ --- an element of the projective plane,

3) $g = \Rt(e_3, l)$ \ is the operator of rotation of $V$ from the axis of 
   $e_3$ to the axis $l$ about the normal to the plane $L(e_3, l)$, \ \
   and we put \ $g = E$ \ for the case of \ $l = e_3$,

4) $g^c$ \ is the change in $\Sym(3)$ by this rotation.

By this construction, $g \in \SO(3)$, so the image of \ $s$ \ resides on the 
orbit of $D$. And we need to prove that \ $s$ \ is a diffeomorphic imbedding 
from the projective plane into $\Sym(3)$, onto the orbit. 

First, $s$ is smooth due to its formula --- the matrix transposition and 
multiplication.\\
Then, notice that \ {\bf 1)} by the point (5) of Lemma Orbit1, each operator 
$A$ on $\Orbit(D)$ has the spectrum \ $\spe = \{1,1,-2\}$,

{\bf 2)} there exists an axis \ $l$ \ for this $A$ which belongs the 
  eigenvalue \ -2. 
  Such an axis is unique, because otherwise, the value \ -2 \ occurs 
  (by Theorem DO) multiple in $\spec(A)$ --- which contradicts the invariant
  value of the multiset \ $\spe$. 

Further, consider the operators $s(l_1), s(l_2)$ for any axes $l_1 \neq l_2$.
These two operators have different axes belonging to the eigenvalue -2.
As it is proved earlier, such an axis is unique for each form in $\Orbit(D)$.
Hence, \ $s(l_1) \neq s(l_2)$, \ which proves \ {\it injectiveness} \ of 
the map $s$.

Let us prove surjectiveness of $s$. \ 
Let $A \in \Orbit(D)$. \ It has the spectrum \ 
$\spe = \spec(A) = \{1,1,-2\}$, \ and there exists a unique axis $l$ for 
$A$ belonging the eigenvalue -2. \ Put \ \ $g = \Rt(e_3,l)$, \ \ and 
let us prove that \ \ $s(l) = g^c D = A$ \ \ by investigating the following 
two cases.

1) $l = e_3$. \  Then, \ $g = E, \ A e_3 = -2 e_3$; \ 
 the 2-dimensional eigenspace \ $V_1 = L(e_1, e_2)$ \ of the eigenvalue
 1 is orthogonal to $e_3$; \ the operator $A$ is identical on the space $V_1$.
 Hence \ $A = D = g^c D = s(l)$. 

2) $l \neq e_3$. \ $g \ e_3 = l$, \ and as $g$ is orthogonal, \ 
$g \ L(e_1, e_2)$ \ is the orthogonal complement \ $l^\bot$ \ to $l$. \ 
The axis $l$ belongs to the eigenvalue -2 for $A$. Hence, by Theorem DO, 
$l^\bot$ is the eigenspace of \ 1 \ for $A$. \ Hence the operator $A$ 
has the matrix \ $\diag(1,1,-2) = D$ \ in the basis \ 
$(g \ e_1, g \ e_2, g \ e_3)$. \ This means \ $g^c D = A$, \ and this 
proves surjectiveness of \ $s$.
\smallskip

\noindent {\bf Remark.} 

1. The idea of parameterizing this orbit by the axis of the eigenvalue $-2$ 
   has been communicated to me by \ S.~Yu.~Orevkov.

2. Even though the projective plane is not a sub-manifold in $V$, 
   the above orbit --- the image of the map $s$ --- is a smooth algebraic 
   variety in $\Sym(3)$.
\smallskip

The orbit also has certain remarkable \ metrical properties. Thus, let us 
show the orbit construction by a circumference family.
\smallskip

\noindent \underbar{Proving (2)} \ \ 
Consider first the 1-orbit of $D$ by the axis $e_1$. 
This is a circumference (residing in its bi-dimensional plane), which we 
denote $\Cir_1$. \ This is visible from the explicit formulae for a change 
in $D$ by a rotation about $e_1$. The general matrix of this 1-orbit is
{\footnotesize
$$\C1G = e_1^\phi D = \pmatrix{1 &  0           &  0\cr
                               0 &  c^2 -2 b^2  &  -3 c b\cr
                               0 &  -3 c b      &  -2 c^2 +b^2},$$
}\noindent 
where \ $c = \cos \phi, \ b = \sin \phi$. \ 
The only diagonal matrices in this circumference are $D$ and \ 
$D_2 = \diag(1,-2,1)$.  
When the angle $\phi$ changes from zero to $\pi/2$, \ the matrix 
$\C1G$ passes (without repetition) the arch from $D$ to $D_2$. 
At this stage, the matrix element $\C1G(2,3)$ first changes
monotonously from zero to $-3/2$, and then changes monotonously back
to zero. When the angle $\phi$ changes from $\pi/2$ to $\pi$, the  
matrix $\C1G$ passes (without repetition) the arch from $D_2$ to $D$; 
this is the ``lower'' half of the 1-orbit, it differs from the upper 
half only in the sign of the element $\C1G(2,3)$.

The subspace $L_1$ of the 1-orbit $\Cir_1$ \ is linearly generated by
the vector \ $\overline{D \ D_2}$ \ and the tangent vector to $\Cir_1$ 
in the point $D$. After norming by the factor $-3$, these two vectors 
produce the basis of the following three matrices for the space $L_1$: 


{\footnotesize
$$M_1 = \diag(0,1,-1), \ \ \ \ M_2 = \pmatrix{0 & 0 & 0\cr
                                              0 & 0 & 1\cr
                                              0 & 1 & 0}.$$
} 
Now, let $l$ be the axis in the space $L(e_1, e_2)$ obtained from $e_1$ 
by a rotation \ $h = \RtA(e_3, \psi)$. \ 
Then the 1-orbit $\Cir(\psi)$ for $D$ by $l$ is a circumference containing 
$D$ and isometric to $\Cir_1$;  this isometry is the operator $h^c$ 
(in the space $\Sym(3)$), and it also 
maps the plane $L_1$ to the plane $L(\psi)$ of the 1-orbit $\Cir(\psi)$. \ 
This property holds due to that the operator $h^c$ conjugates the 
subgroups $G(e_1)$ and $G(l)$. Therefore, for these two 1-orbits it holds
the equation
$$           \Orbit(l) D = \ h^c \Orbit(e_1) (h^{-c} D), 
$$
where \ $\Orbit(l) A$ \ denotes the orbit of $A$ by the changes by rotations
about the axis $l$. 
But in this case, $h$ is a rotation about $e_3$, and $D = \diag(1,1,-2)$. \ 
Hence $h$ commutes with the operator $D$, and the relation between these 
two 1-orbits is simplified:  
                            $$\Orbit(l) D = \ h^c \Orbit(e_1) D.$$ 
\noindent
That is these 1-orbits are conjugated by $h$. 

So, it appears that the orbit for $D$ contains the union of a smooth
single-parameter family of circumferences having a common point $D$ and 
obtained from $\Cir_1$ by various rotations about $D$.\\  
{\it The point (2) is proved}. 

Let us also note (informally, for intuition) that the plane $\Pi(\psi)$ of 
1-orbit in this family turns dependently on $\psi$ in various directions in 
a five-dimensional space.
\smallskip

\noindent \underbar{Proving (3)} \ \ 
We need to prove that the orbit of the above matrix $D$ is not contained 
in any four-dimensional plane. 
The simplest way to see this is to consider any embracing plane for the 
orbit of \ $D = \diag(0,0,1)$. \ 
Because by the statement (7) of Lemma Orbit1, the latter orbit differs 
metrically from the former only by some scalar shift and homothety by some
non-zero factor. 

It is easy to compute that the tangent space to the orbit in the point $D$ 
is generated by the vectors 
{\footnotesize
$$T_1 = \pmatrix{0 & 0 & 0\cr
                 0 & 0 & 1\cr 
                 0 & 1 & 0}, \ \ \ 
  T_2 = \pmatrix{0 & 0 & 1\cr
                 0 & 0 & 0\cr
                 1 & 0 & 0}.
$$}\noindent
These are the results of changes in $D$ by infinitesimal rotation about 
$e_1$ and about $e_2$ respectively, they are obtained by commuting the 
matrix $D$ with two (anti-symmetric) matrices from a basis of the Lie 
algebra $\so(3)$.

Also the orbit contains the matrices  \ $D_2 = \diag(0,1,0)$ \ and \ 
                                        $D_3 = \diag(1,0,0)$. 
\ \ 
Respectively, the tangent space in $D_3$ contains the vector
{\footnotesize
$$B = \pmatrix{0 & 1 & 0\cr
               1 & 0 & 0\cr
               0 & 0 & 0}.$$
}
Denote \ \ $D_2\prime = D_2 - D, \ \ D_3\prime = D_3 - D$. 

If a plane $L$ contains the orbit for $D$, then $L$ contains the
point $D$ and also the directions \ 
 $T_1$, \ $T_2$, \ $B$, $D_2\prime = \diag(0,1,-1)$, \ 
 $D_3\prime = \diag(1,0,-1)$. \ \ 
The vectors \ $T_1$, $T_2$ and $B$ \ are mutually orthogonal (recall the 
s-metric).
The vector $D_i\prime$ is orthogonal to $T_1$, $T_2$ and $B$ \ for 
$i = 2, 3$, \ 
and the vectors \ $D_2\prime$ and $D_3\prime$ \ constitute a staircase
matrix of rank 2. So: the above five directions are linearly 
independent and belong to each plane containing $\Orbit(D)$.\\
{\bf The lemma is proved}.

\subsection{Proof for the structure of $\cal{M}$ for $n = 3$}
\label{sec-DSForm:structure3}

By statement (8) of Lemma Orbit1, $\cal{M}$ is the straight cylinder over 
the surface restriction $\cal{M}$$_0$. It remains to describe the 
geometric structure of $\cal{M}$$_0$. Each matrix in $\cal{M}$$_0$ has 
spectrum of the kind $\{\lambda, \lambda, -2 \lambda\}$, because the three 
eigenvalues satisfy the conditions of zero sum and existence of a multiple 
value. Therefore, by Theorem DO, the diagonal matrices in $\cal{M}$$_0$ 
are all the matrices of the kind \ 
$D(\lambda)   = \diag(\lambda, \lambda, -2 \lambda)$ \ or  
$D_2(\lambda) = \diag(\lambda, -2 \lambda, \lambda)$ \ or \ 
$D_3(\lambda) = \diag(-2 \lambda, \lambda, \lambda)$.
\\
{\it The statement (4) of Theorem DS3 is proved}.
\smallskip

\noindent {\bf Remark.} \ \ 
{\bf 1)} For each $\lambda$, the matrices $D(\lambda), D_2(\lambda)$ and  
   $D_3(\lambda)$ \ belong to the same orbit 
   (by the statement 5 of Lemma Orbit1). 

{\bf 2)} The points $D(1)$ and $D(-1)$ \ are central-opposite to each other 
         on the sphere $S_4$ (to which $D(1)$ belongs), 
         but they belong to different orbits.

  The conic set $\cal{M}$$_0$ is union of the ``positive'' half-cone ---
  union of orbits of the family $D(\lambda)$ for $\lambda \geq 0$, \ 
  and the ``negative'' half-cone --- union of orbits of the family 
  $D(\lambda)$ for $\lambda \leq 0$.

{\bf 3)} Intersection of \ $\cal{M}$$_0$ \ with the sphere \ $S_4$ \  
         (to which $D(1)$ belongs), \ is\\ $\Orbit(D(1)) \cup \Orbit(D(-1))$.
\smallskip

For any form $A \in \cal{M}$$_0$$ \ \ \spec(A)$ is proportional to the 
multiset $\{1,1,-2\}$ with a factor $\lambda$ depending on $A$. Therefore, \ 
$\Orbit(A) = \Orbit(\lambda(A) D) = \lambda(A) \cdot \Orbit(D)$.\linebreak
Here \ $D = \diag(1,1,-2)$, \ $\lambda(A)$ are all real numbers, \ 
and the last equality is true due to that \ 
$g^c \lambda D = \lambda g^c D$ \ for all invertible $g \in {\cal L}(n)$. 

Therefore, $\cal{M}$$_0$  is the cone with its vertex in zero and its base 
being $\Orbit(D)$.\\ 
And by Lemma ST of Section \ref{sec-DSForm:orbit3},
$\Orbit(D)$ is a diffeomorphic image of the projective plane, and it resides 
on the sphere $S_4$ with center in zero.

\newpage
\noindent {\bf Theorem.} \ \ 
         Zero matrix is the only singular point on the surface $\cal{M}$$_0$.

\noindent {\bf Proof.}\\ 
Any non-zero form $A \in \cal{M}$$_0$ belongs to $\cal{MH}$. By \cite{Ar1} 
(see Section \ref{sec-Ar} of this paper), such $A$ is a regular point on 
$\cal{M}$. \ 
As $\cal{M}$ is a straight cylinder over $\cal{M}$$_0$, \ $A$ is also 
regular on $\cal{M}$$_0$.

Why zero is singular? In our case, the cone base is a 2-dimensional surface 
injected in a complex way into a 4-dimensional sphere. So, we cannot see 
singularity of the vertex by simply referring to the picture of a straight 
circular cone. Let us prove singularity by finding the enough number of 
linearly independent tangent vectors.

The vertex ($0$) is singular because \  
1) $\cal{M}$$_0$ is a conic set with vertex in zero,\\ 
2) we present certain four forms $A_i$ on the orbit of $D = \diag(1,1,-2)$ 
such that their axes $\overline{0 A_i}$ are linearly independent (and
reside in $\cal{M}$$_0$).
   
The two of these forms are \ $D = \diag(1,1,-2)$ \ and \ 
$D_2 = \diag(1,-2,1)$. \ 
The orbit of $D$ contains the 1-orbit of $D$ by the axis $e_1$. 
This is a circumference containing the matrices $D$ and $D_2$, so that $D$ 
corresponds to the rotation angle $0$, \ and $D_2$ corresponds the angle 
$\pi/2$. \ 
We also need a couple of non-diagonal matrices on the orbit. 
Thus, \ $e_1^{\pi/4} D$ \ is the matrix \ 
{\footnotesize
$$M_1 = \pmatrix{1 & 0 & 0  \cr
                 0 & a & b  \cr
                 0 & b & -a-1}, \ \ \ a = -1/2, \ b = -3/2$$
}\noindent
(recall that \ $e_i^\phi$ \ is the change by the rotation operator 
$\RtA(e_i, \phi)$).

Similarly, the 1-orbit for $D$ by $e_2$ contains the matrix
{\footnotesize
$$M_2 = \pmatrix{a  &  0  &  b \cr
                 0  &  1  &  0 \cr  
                 b  &  0  &  -a-1}$$ 
}\noindent 
--- with the same \ $a$ and $b$. \
To verify linear independence of the symmetric matrices \ 
$D$, $D_2$, $M_1$, $M_2$, \ represent each of them in the 
row--vector form: \ $[X(i,j) | \ 1 \leq i \leq j \leq 3]$. \ 
The sub-diagonal part of each matrix is skipped in this representation, 
because it equals the super-diagonal part.
This produces one matrix of the four rows:   
{\footnotesize
$$\pmatrix{1 & 0 & 0 & 1 & 0  & -2   \cr 
           1 & 0 & 0 & 0 & -2 & 1    \cr
           1 & 0 & 0 & a & b  & -a-1 \cr
           a & 0 & b & 1 & 0  & -a-1}.$$
}
The first loop of clearing of the first column by the staircase form 
algorithm produces the matrix
{\footnotesize
$$\pmatrix{1 & 0 & 0 & 1   &  0 & -2   \cr
           0 & 0 & 0 & -1  & -2 & 3    \cr
           0 & 0 & 0 & a-1 & b  & -a+1 \cr
           0 & 0 & b & 1-a & 0  & a-1}.$$
}
\noindent
Further, the fourth row is moved to the second place, and to the third
row it is added the second row multiplied by \ $a-1$. This produces 
the matrix
{\footnotesize
$$\pmatrix{1  &  0  &  0  &  1   &  0   &  -2  \cr
           0  &  0  &  b  &  1-a &  0   &  a-1 \cr
           0  &  0  &  0  &  -1  &  -2  &  3   \cr
           0  &  0  &  0  &  0   &  b_2 &  a_2}.$$
}\noindent 
Here \ \ $b \neq 0, \ \ \ b_2 = b - 2(a-1) = -3/2 - 2(-3/2) \neq 0$. \ \ \ 
Hence this is a staircase matrix of rank 4, and hence the matrices
$D$, $D_2$, $M_1$, $M_2$ \ are linearly independent, and their
corresponding axes are linearly independent.
So, we present the four linearly independent axes on the surface 
$\cal{M}$$_0$ which intersect in zero. Therefore, zero is a singular 
point in $\cal{M}$$_0$. \ 
Let us explain --- why. It is proved earlier that all other points 
of this surface are regular. So, the tangent space to each non-zero point 
on $\cal{M}$$_0$ has dimension 3. To prove by contradiction, suppose that 
zero in regular on $\cal{M}$$_0$. 
Then, by the definition of a regular point, some neighbourhood of zero in 
$\cal{M}$$_0$ is diffeomorphic to a ball in $\mathbb{R}^3$. Due to its 
dimension, such a neighbourhood cannot contain the above four linearly 
independent axes coming out from zero.
This contradiction shows that zero is a singular point in 
$\cal{M}$$_0$.\\
{\it The theorem is proved.}
\smallskip

So, 0 is the only singular point on $\cal{M}$$_0$. By Theorem DS, \ $\cal{M}$ 
is the straight cylinder over $\cal{M}$$_0$ with the line $\Scal$ as element. 
This proves the last remaining statement in Theorem DS3: 

 singular points on $\cal{M}$ are the scalar matrices and only them.
\\
{\bf Theorem DS3 is proved.}

\newpage
\section{Appendix}
\label{sec-attach}

The proofs in Appendix do \ not \ use the statements being proved in the 
main part of the paper.

All the statements presented in Appendix are evident or known, 
their proofs are given for a formal reason.
\smallskip

\noindent {\bf Lemma (Stab).} \ \ 
Let \ $M \in \cal{MH}$. 

1) The set $V_1$ of the vectors belonging to the multiple eigenvalue
   (which is unique) for $M$ is a subspace of dimension 2.
 
2) The connected component $G_e$ of unity in the stabilizer $\St(M)$ 
   in $\SO(n)$ \ is a single-parameter subgroup of the operators which 
   restriction on $V_1$ are rotations of this 2-dimensional space about 0, 
   and which restriction to the orthogonal complement to $V_1$ is identity.

\noindent {\bf Proof.}\\ 
Let $\lambda$ be the (unique) multiple eigenvalue (of the multiplicity 2) of 
the operator $M$. By Theorem DO, there exists a basis \ $B$ \ in the space 
$V$ \ consisting of the eigenvectors of $M$, and there are exactly two 
vectors in $B$ belonging to the eigenvalue $\lambda$. 
Denote these vectors \ $u$ and $v$; \ \ 
denote \ $B_2 = B \backslash \{u,v\}$, \ \ $V_1 = L(u,v)$, \ \ 
$V_2 = L(B_2)$.

The condition \ $g \in \St(M)$ \ means that $g$ commutes with $M$: \ 
$g M = M g$. \ 
Therefore, if a vector $v$ belongs to any eigenvalue $\mu$ for $M$, then 
it holds \ \ $M g v = g M v = g \mu v = \mu g v$. \ \ 
Hence, each eigenspace for $M$ is mapped by $g$ on itself. 

Let \ $g_1$ and $g_2$ \ be the restriction of \ $g$ \ to $V_1$ and to 
$V_2$ \ respectively.

Each vector in \ $B_2$ \ belongs to a 1-dimensional eigenspace for $M$, 
hence such a vector is an eigenvector for \ $g$. 
Further, an operator from $\SO(n)$ can have only the eigenvalue 1 or -1.
Therefore, in the basis $V_2$ the operator $g_2$ has a diagonal matrix 
with each diagonal value being 1 or -1. 

Further, $g$ maps the subspace $V_1$ onto itself, and it is an orthogonal
operator ($g_1$) on $V_1$. 
The operator $M$ is scalar on $V_1$, hence its restriction to $V_1$ 
commutes with each linear operator on $V_1$.

$\det(g) = 1$. \  Hence, if the multiplicity \ $m$ \ of \ -1 in $\spec(g_2)$
is even, then $g_1$ is oriented positively, and $g_1$ is a rotation about 
0 at some angle. \ If $m$ is odd, then $g_1$ is oriented negatively, and 
$g_1$ is a rotation about 0 composed with some reflection mapping.

$\St(M)$ \ is an algebraic Lie group, as it is widely known 
(our reference for the Russian language is \cite{VinO}, head I, Theorem 1, 
and the following exercise).
It consists of several connected components.
Let $H$ by the set of operators from $\St(M)$ which are identical on $V_2$.
Then \ $H \subseteq G_e$, \ because it is proved that for each operator 
from $H$ its restriction on $V_1$ is a rotation of a 2-dimensional plane 
about 0, \ and the set of such operators contains unity, is connected and 
smoothly parameterized by the angle of rotation.

It is proved above that each operator $g$ in $St(M)$ has on $V_2$ in the 
basis $B_2$ a diagonal matrix $D(g)$ with each diagonal value being 1 or 
-1. The set of such matrices $D(g)$ is finite and discrete. 
Therefore, for each operator $g$ in the connected component $G_e$ the 
operator $g_2$ has the only eigenvalue 1, hence $g_2$ is identical. 
Therefore \ $H = G_e$ \ --- the connected component of unity.
\\
{\bf The lemma is proved.}
\medskip

\noindent {\bf Lemma (AlgSO)} (known). \  

{\bf (1)} The tangent space $T$ to the surface $\SO(n)$ in the point $E$
    (the Lie algebra $\so(n)$) consists of all the antisymmetric matrices.

{\bf (2)} The tangent space $T$ to the orbit of a diagonal matrix $D$ by 
  the changes of $\SO(n)$ in the point $D$ \ consists of the commutators of 
  $D$ with all the antisymmetric matrices.

\noindent {\bf Proof.} \ 
Denote \ $G = \SO(n)$. \

\noindent \underbar{Proving (1)} \ \
By definition, $so(n)$ consists of the derivatives  \ $g^\prime(0)$ \ at zero 
of all smooth curves $g(t)$ on the surface $G$, such that \ $g(0) = E$. \
By the property of the transposed orthogonal matrix, it holds \ \
$g(t) g(t)^* = g g^{-1} = E$. \ \ Differentiating this equation at $t = 0$, 
taking in account that transposing a matrix is a linear map, and by the 
rule of derivative of a product of linear operators depending on $t$, we 
obtain: 
$$0 = \ (g(t) g(t)^*)^\prime(0) = \                                       
              g^\prime(0) g(0)^* + g(0) g^\prime(0)^* \ = \                     
              g^\prime(0) + g^\prime(0)^* = \ 0.
$$
\noindent
This proves the statement (1) of that an arbitrary matrix ($g^\prime(0)$) 
in $so(n)$ is antisymmetric.

\noindent \underbar{Proving (2)} \ \
By definition, $T$ consists of the derivatives at $t = 0$ of all the curves 
of the kind \ \ $\gamma(t) = g(t)^c D$, \ \ where $g(t)$ is a smooth curve 
on $G$ such that \ $g(0) = E$. \ By the rule of derivative of a product of 
linear operators depending on $t$, we have

   $(d/dt \ \gamma)(0) = \ d/dt \ (g(t) \ D \ g(t)^*)(0) = \                   
    g^\prime(0) \ D \ g(0)^* \ + \ g(0) \ D \ g^\prime(0)^* =$

   $g^\prime(0) D \ + \ D \ g^\prime(0)^*$    \hspace{100mm}  (1).

Differentiation at zero of the true equation \ \ $g(t) g(t)^* = E$ \ \ 
(for an orthogonal $g$) gives the equation \ \ $g^\prime(0)^* = -g^\prime(0)$. 
\ \
Therefore, the equation (1) is equivalent to \ 
$(d/dt \ \gamma)(0) = [g^\prime(0), \ D]$ \ --- a commutator of two matrices.
In our case, $g^\prime(0)$ is an arbitrary element of the Lie algebra 
$\so(n)$. By the statement (1), this algebra consists of all the 
antisymmetric matrices.\\ 
{\bf The lemma is proved.}
\medskip

\noindent {\bf Lemma A0} (probably, known). \ 

For a diagonal quadratic matrix $D$ and an anti-symmetric matrix $A$ of 
the same size, the commutator $[D, A]$ is a symmetric matrix having zero 
diagonal.

\noindent {\bf Proof.}\\ 
Symmetry of the result matrix is proved by the following equations:
$$ (D A - A D)^* =  (D A)^* - (A D)^* =  A^* D^* - D^* A^* =  
   (-A) D - D (-A) =  D A - A D.
$$
Further, it is sufficient to prove that \ \ \ 
{\bf (a)} \ $D A$ \ \ has zero diagonal \ and\\
{\bf (b)} \ $A D$ \ has zero diagonal. \ \
Also (a) implies (b), because \ $A D = - (D A)^*$, \
and the maps of transposition and negation $-M$ for any matrix $M$ preserve 
the property of having zero diagonal.
By the matrix product definition, for each position $i$ on the diagonal it
holds

\centerline{$(D A)_{i,i} = \ \sum_{j=1}^n \ D_{i,j} A_{j,i} = \                
                                            D_{i,i} A_{i,i} = 0$.
}\noindent
This is because $D_{i,j} = 0$ \ for $i \neq j$, \ and $A_{i,i} = 0$.\\
{\bf The lemma is proved.}

\noindent {\bf Illustration:} \ for 
{\footnotesize
$$D = \diag(a, b), \ \ A = \pmatrix{0  & 1 \cr
                                    -1 & 0}
$$}\noindent it holds \ $[A, D] =$
{\footnotesize
$$\pmatrix{0  & 1 \cr
           -1 & 0} \times
  \pmatrix{a & 0 \cr
           0 & b}
  \ - \
  \pmatrix{a & 0 \cr
           0 & b}     \times
  \pmatrix{0  & 1 \cr
           -1 & 0}  \ = \
  \pmatrix{0  & b \cr
           -a & 0}  -
  \pmatrix{0  & a \cr
           -b & 0} \ = \
  \pmatrix{0   & b-a \cr
           b-a & 0}.
$$}
\smallskip

\noindent {\bf Lemma A1} (known).
 
For any polynomial map \ $F : \ R^n \rightarrow R^m$, \ the counter-image
$\cal{M}^\prime$ of an algebraic set $\cal{M}$ is an algebraic set.

\noindent {\bf Proof.} \ 
Let $\cal{M}$ be the set of zeroes for a polynomial system 
$\{p_1, \ldots, p_k\}$. \ 
Then $\cal{M}^\prime$ is the set of zeroes for the set of compositions 
$\{p_1(F(X)), \ldots, p_k(F(X))\}$. \ Indeed, if $A$ belongs to 
$\cal{M}^\prime$, then $F(A)$ belongs to $\cal{M}$; then \ $p_i(F(A)) = 0$ \ 
for all $i$.
 
Conversely: if \ $p_i(F(A)) = 0$ \ for all $i$,  then $F(A)$ belongs to
$\cal{M}$, hence $A$ belongs to $\cal{M}^\prime$ \ 
(by the definition of counter-image of a set).
\medskip

\noindent {\bf Lemma A2} (known). 
  
For a surjective polynomial map from an algebraic set $\cal{M}$ onto an
algebraic set $\cal{M}^\prime$, \
if $\cal{M}$ is irreducible, then $\cal{M}^\prime$ is irreducible.

\noindent {\bf Proof.}\\
To prove by contradiction, suppose that $\cal{M}^\prime =$ $M_1 \cup M_2$ \ 
for algebraic sets $M_1$ and $M_2$, each different from $\cal{M}^\prime$.
By Lemma A1, \ the counter-images $C_1$ and $C_2$ of $M_1$ and of $M_2$ \ 
respectively are algebraic sets. Also it holds \ $\cal{M} =$ $C_1 \cup C_2$. \ 
Further, $C_1 \neq \cal{M}$, because $M_1 \neq \cal{M}^\prime$. \ 
For a similar reason, \ $C_2 \neq \cal{M}$. \ So, this representation \ 
$\cal{M} =$ $C_1 \cup C_2$ \ contradicts irreducibility of $\cal{M}$.
\\
{\bf The lemma is proved.}


\newpage
\begin{thebibliography}{Il1}
%
\bibitem[Ar1]{Ar1} V.~I.~Arnold, K.~Vogtmann, A.~Weinstein. \ 
{\it Mathematical Methods of Classical Mechanics}. Second edition.
Springer Science + Business, Media Inc. 1978, 1989.\\
Translation of ``Matematicheskie metody classicheskoi mekhaniki'', 1974.

\bibitem[Ar2]{Ar2} V.~I.~Arnold. \ \ 
{\it Relatives of the quotient of the complex projective plane by 
     the complex conjugation}. \ 
English and Russian versions are available as files: \ 
\verb#http://www.pdmi.ras.ru/~arnsem/Arnold/arn-papers.html# (1998).\\
Russian version is also published in Proceedings of the Steklov's
Mathematical Institute, 1999, volume 224, pages 56 -- 67.

\bibitem[Bo]{Bo} A.~Borel. {\it Linear Algebraic Groups}.\\ 
Published by W.~A.~Benjamin, New York--Amsterdam, 1969.

\bibitem[Il]{Il} N.~V.~Ilyushechkin. \ 
{\it Discriminant of the Characteristic Polynomial of a Normal Matrix.}
Mathematical Notes, 1992, {\bf 51}:3, \ Russian version: pages 16--23, \ 
English version: pages 230--235.

\bibitem[Il2]{Il2} N.~V.~Ilyushechkin. \ 
 {\it On the relations between the summands of a symmetric matrix 
      discriminant}. \ \ 
 To be published in ``Mathematical Notes''.

\bibitem[Me]{Me} S.~D.~Mechveliani. \ 
{\it Geometry of the Variety of Quadratic Forms with Multiple Eigenvalues}. \ 
Manuscript in the electronic library of Cornell University, 2009, 2010. \ \ 
\verb#http://arxiv.org/abs/0907.3293#, \ .pdf, .ps files. 

\bibitem[Sha]{Sha} I.~R.~Shafarevich. \ {\it Basic Algebraic Geometry: 1.}\\
Springer-Verlag, Berlin, Heidelberg, 1994.

\bibitem[VinO]{VinO} E.~B.~Vinberg, A.~L.~Onishik. \ 
{\it Seminar on Lie Groups and Algebraic Groups}. \ 
(``Seminar po Gruppam Li i Algebraicheskim Gruppam''). In Russian.
``Nauka'', Moscow, 1988.

\bibitem[Wa]{Wa} B.~L.~Van Der Waerden. \ {\it Algebra I}, {\it Algebra II}. 
Springer-Verlag, Berlin Heidelberg, New York, 1971, 1967. 
\end{thebibliography}
\end{document}